\documentclass[a4paper,12pt,leqno]{article}

\usepackage{amsmath}
\usepackage{amssymb}
\usepackage{amscd}
\usepackage{graphicx}

\newcommand{\vect}[1]{{\textnormal{\mathversion{bold}$#1$}}}
\newtheorem{theorem}{Theorem}

\newcommand{\qed}{\hfill $\square$  \vspace{24pt}}

\begin{document}

\title{On the Scarf-Hirota model\\ in the price-scaled price adjustment process}
\author{Tatsuro Yamamoto,  Masanari Asano \thanks{\texttt{e-mail address : asano@is.noda.tus.ac.jp}} , \\ Yoshio Togawa and Masanori Ohya \\ \\ Dept.\ of Information Sciences, 
  Faculty of Science \& Technology,\\
  Tokyo University of Science,
  Noda City, Chiba 278, Japan.}
\date{}
\maketitle
\begin{abstract}
Hirota's results given in (Hirota.M.,1981) on the asymptotically stability are generalized to the price-scaled price adjustment process.
\end{abstract}


\newpage

\section{Introduction}

We first review the Scarf-Hirota model for exchange economy. The Scarf-Hirota model is an exchanged economy which consists of three consumers and three commodities. The consumers utility functions are 
\begin{eqnarray*}
u_1(\vect{x_1}) &=& \min\{x_{21},x_{31}\}\\
u_2(\vect{x_2}) &=& \min\{x_{12},x_{32}\}\\
u_3(\vect{x_3}) &=& \min\{x_{13},x_{23}\}
\end{eqnarray*}
where the vector $\vect{x_i}=(x_{1i},x_{2i},x_{3i})^T$ $(i=1,2,3)$ denotes the $i$-th consumer\textbf{'s} consumption function ("$T$" is the transpose of the vector).

The total excess demand functions are given by
\begin{eqnarray*}
E_1(\vect{p})&=&\frac{\vect{p}^T\vect{a}_2}{p_3+p_1}+\frac{\vect{p}^T\vect{a}_3}{p_1+p_2} - (a_{11}+a_{12}+a_{13})\\
E_2(\vect{p})&=&\frac{\vect{p}^T\vect{a}_3}{p_1+p_2}+\frac{\vect{p}^T\vect{a}_1}{p_2+p_3} - (a_{21}+a_{22}+a_{23})\\
E_3(\vect{p})&=&\frac{\vect{p}^T\vect{a}_1}{p_2+p_3}+\frac{\vect{p}^T\vect{a}_2}{p_3+p_1} - (a_{31}+a_{32}+a_{33}).
\end{eqnarray*}
where $\vect{a_i}$ $(i=1,2,3)$ denotes the initial endowment vector of $i$-th consumer.

In the Scarf-Hirota model, the matrix
\[
A=[\vect{a_1},\vect{a_2},\vect{a_3}].
\]
satisfies the following condition: For a vector  $\vect{u}=[1,1,1]^T$, 
\begin{equation}
\mbox{``{\em Condition A}:}\ \ \ \ \ A\vect{u}= \vect{u},\quad \vect{u}^TA=\vect{u}^T\mbox{"}.
\end{equation}

The \textit{t\^{a}tonnement} process in the Scarf-Hirota model is described by the differential equation
\begin{equation}
\frac{dp_i}{dt} = E_i(\vect{p}), \ i=1,2,3, \label{classical}
\end{equation}
which has a fixed point $\vect{p^*}$.

Hirota found a condition of asymptotic global stability of the fixed point $\vect{p^*}$, where the condition is described by an inequality of elements of matrix $A$.

In this paper, the \textit{T\^{a}tonnement} process we consider is the price-scaled price adjustment process (Anderson C.M., Plott C.R., Shimomura K.I., Granat S. 2004)
\begin{equation}
\frac{dp_i}{dt} = p_iE_i(\vect{p}),\ i=1,2,3, \label{price scale}
\end{equation}
and more generally
\begin{equation}
\frac{dp_i}{dt} = p_i^\gamma E_i(\vect{p}),\ i=1,2,3, \label{lambda}
\end{equation}
where $\lambda$ is an arbitrary nonnegative number
We will study the local and global asymptotically stability for the above differential equations.

In Section 2, we prepare the notations for the Scarf-Hirota model and find a first integral in general form. The global stability is discussed in Section 3, and the local stability is considered in Section 4. The orbits close to the boundary are examined in Section 5.

\section{Preliminary }

The following notations are used:

  (ⅰ) Initial endowment matrix is described by $A=[\vect{a}_1,\vect{a}_2,\vect{a}_3]^T$.

  (ⅱ) Nominal price vector is described by $\vect{p}=[p_1,p_2,p_3]^T$.

  (ⅲ) Equilibrium price vector is described by $\vect{p*}=[1,1,1]^T$.

  (ⅳ) $B=[\vect{b}_1,\vect{b}_2,\vect{b}_3]=\left[
\begin{array}{lll}
0&1&1\\
1&0&1\\
1&1&0\\
\end{array}\right]$,

where $\vect{b}_1=[0,1,1]^T,\quad \vect{b}_2=[1,0,1]^T,\quad \vect{b}_3=[1,1,0]^T$.

  (ⅴ) $\vect{f}=[f_1,f_2,f_3]^T$, where 
  \[ f_h=\frac{\vect{p}^T \vect{a}_h}{\vect{p}^T\vect{b}_h} \].

The maximum utility of the consumer $h=(1,2,3)$ is achieved by
\[
\vect{x}_h=f_h\vect{b}_h.
\]
Then the total excess demand function is given by
\begin{equation*}
\vect{E}(\vect{p})=B\vect{f}-A\vect{u}.
\end{equation*}
where $\vect{u}=[1,1,1]^T$.
One can write the matrix $A$ satisfying the condition (A) as
\begin{equation*}
A=\left[
\begin{array}{ccc}
a_{11} & a_{12} & 1-a_{11}-a_{12}\\
a_{21} & a_{22} & 1-a_{21}-a_{22}\\
1-a_{11}-a_{21} & 1-a_{12}-a_{22} & a_{11}+a_{12}+a_{21}+a_{22}-1\\
\end{array}\right]
\end{equation*}
(Hirota~\cite{Hirota}).
We rewrite the above matrix $A$ as
\begin{equation*}
A=\left[
\begin{array}{ccc}
d_1 & d_3+L & d_2+K\\
d_3+K & d_2 & d_1+L\\
d_2+L & d_1+K & d_3
\end{array}\right],
\end{equation*}
with the parameters $d_1\geq 0,d_2\geq 0,d_3\geq 0,K,L$ satisfying
\begin{equation}\label{dKL}
d_1+d_2+d_3+K+L=1.
\end{equation}
According to the condition (A), these parameters satisfy the following conditions:
\begin{equation*}
a_{11}=d_1,\quad a_{22}=d_2,\quad a_{33}= d_3,\quad a_{12}=d_3 + L,
\end{equation*}
\begin{eqnarray*}
a_{13}&=&1-a_{11}-a_{12}=d_2+K,\\
a_{32}&=&1-a_{12}-a_{22}=d_1+K,\\
a_{23}&=&1-a_{13}-a_{33}=d_1+L,\\
a_{21}&=&1-a_{22}-a_{23}=d_3+K,\\
a_{31}&=&1-a_{11}-a_{12}=d_2+L.
\end{eqnarray*}
In Section 3, it will be shown that the condition of asymptotic global stability is given as
\begin{equation}
H \equiv 4(d_1d_2+d_2d_3+d_3d_1) - (d_1+d_2+d_3)^2> 0
\end{equation}
In Section 4, the condition of the asymptotic local stability is given by
\begin{equation}
\hat{H} \equiv H+(K-L)^2 > 0
\end{equation}
The main purpose of this paper is to discuss the price-scaled price adjustment process (\ref{price scale}).
For such a discussion, we use more general form (\ref{lambda})
where $\gamma$ is a nonnegative number. Note that $\gamma=0$ corresponds to the classical case (\ref{classical}) and $\gamma=1$ corresponds to the price-scaled price adjustment process (\ref{price scale}).

In case of $\gamma=0$, differential equations systems (\ref{classical}) has first integral $p_1^2+p_2^2+p_3^2$.

This integral is generalized for arbitrary$\gamma>0$
\begin{theorem}
In the Scarf-Hirota model
\begin{equation*}
g_\gamma = \left\{
\begin{array}{cl}
p_1 ^{2-\gamma} + p_2 ^{2-\gamma} + p_3^{2-\gamma}, & \mbox{if } \gamma \ne 2\\
p_1p_2p_3, & \mbox{if }\gamma = 2
\end{array}
\right.
\end{equation*}
is the first integral of the ordinary differential equations (\ref{lambda}).
\end{theorem}
\textbf{Proof}: For $\gamma\ne 2$, one has
\[
\frac{dg_\gamma}{dt}=\sum_{i=1}^{3} \frac{\partial g_\gamma}{\partial p_i}\frac{dp_i}{dt} = \sum_{i=1}^3  (2-\gamma)p_i^{1-\gamma} p_i^\gamma E_i = (2-\gamma)\sum_{i=1}^3 p_i E_i
\],
which is equal to 0 due to the Walras law.\\
For $\gamma=2$, one has
\[
\frac{dg_\gamma}{dt}=\sum_{i=1}^{3} \frac{\partial g_\gamma}{\partial p_i}\frac{dp_i}{dt} = p_1p_2p_3\sum_{i=1}^3  p_i^{-1} p_i^2 E_i = p_1p_2p_3\sum_{i=1}^3 p_i E_i = 0.\ \mbox{\qed}
\]

In the exchange economy of n commodities, the first integral is
\begin{equation*}
g_\gamma = \left\{
\begin{array}{cl}
p_1 ^{2-\gamma} + \cdots + p_n^{2-\gamma}, & \mbox{if } \gamma \ne 2\\
p_1p_2\cdots p_n, & \mbox{if }\gamma = 2
\end{array}
\right.
\end{equation*}
The solution is constrained on a $(n-1)$-dimensional surface $g_\gamma=C$

\section{Global Stability}

In this section, we discuss the global stability of the differential equations system (4). Hirota~\cite{Hirota} studied the case of $\gamma =0$ by finding Lyapunov function. Let us introduce the following function of arguments $p_1$, $p_2$ and $p_3$; 
\[
\phi_\gamma = \left\{
\begin{array}{ll}
p_1^{-\gamma} + p_2^{-\gamma} + p_3^{-\gamma} & \mbox{if } \gamma > 0\\
-p_1p_2p_3 & \mbox{if } \gamma = 0
\end{array}
\right.
\]
which is the Lyapunov function including the one that Hirota defined. If the vector {\textnormal{\mathversion{bold}$p$}} satisfies the differential equations (4),
\[
\frac{dp_i}{dt} = p_i^\gamma E_i, \quad i=1,2,3\ ,
\]
the time differentiation of $\phi_\gamma$ is
\begin{eqnarray*}
\frac{d\phi_\gamma}{dt} &=& 
\frac{d\phi_\gamma}{dp_1}p_1^\gamma E_1 +
\frac{d\phi_\gamma}{dp_2}p_2^\gamma E_2 +
\frac{d\phi_\gamma}{dp_2}p_2^\gamma E_2\\
&=& \left\{
\begin{array}{rl}
-\frac{\gamma}{p_1p_2p_3}\left(p_2p_3E_1 + p_1p_3E_2 + p_1p_3E_3\right) & \mbox{if } \gamma > 0\\
-\left(p_2p_3E_1 + p_1p_3E_2 + p_1p_3E_3\right) & \mbox{if } \gamma = 0.
\end{array}\right.
\end{eqnarray*}
Therefore, the value of $\phi_\gamma$ decreases along the orbit of {\textnormal{\mathversion{bold}$p$}} when the inequality  
\[
p_2p_3E_1+p_1p_3E_2+p_1p_3E_3 ( =  [p_2p_3,p_1p_3,p_1p_3] {\textnormal{\mathversion{bold}$E$}} ) > 0
\]
is satisfied. ( Here $p_1>0,\ p_2>0$ and $p_3>0$ are assumed.)\\
\textbf{Lemma 1}
\[
[p_2p_3,p_1p_3,p_1p_2] {\textnormal{\mathversion{bold}$E$}}
\]
is expressed in the quadratic form
\[
[p_2p_3,p_1p_3,p_1p_2] {\textnormal{\mathversion{bold}$E$}} = {\textnormal{\mathversion{bold}$p$}}^{T} (A- \frac{1}{2} B) {\textnormal{\mathversion{bold}$p$}},
\]
where $A$ and $B$ are the matrices given in the previous section.

\noindent
\textbf{Proof}: Note ${\textnormal{\mathversion{bold}$E$}}=B{\textnormal{\mathversion{bold}$f$}}-A{\textnormal{\mathversion{bold}$u$}}$. We have \\
\begin{eqnarray*}
[p_2p_3,p_1p_3,p_1p_2]B {\textnormal{\mathversion{bold}$f$}} &=& p_1(\vect{p}^T{\textnormal{\mathversion{bold}$a_1$}}) +
p_2({\textnormal{\mathversion{bold}$p$}}^T {\textnormal{\mathversion{bold}$a_2$}}) +
p_3({\textnormal{\mathversion{bold}$p$}}^T {\textnormal{\mathversion{bold}$a_3$}})\\
&=& {\textnormal{\mathversion{bold}$p$}}^T A {\textnormal{\mathversion{bold}$p$}},
\end{eqnarray*}
because of 
\begin{eqnarray*}
[p_2p_3,p_1p_3,p_1p_2]B &=& 
[p_1p_3+p_1p_2,p_2p_3+p_1p_2, p_2p_3+p_1p_3]\\
&=&
(p_1({\textnormal{\mathversion{bold}$p$}}^{T}{\textnormal{\mathversion{bold}$b_{1}$}}),p_2({\textnormal{\mathversion{bold}$p$}}^T {\textnormal{\mathversion{bold}$b_2$}}),p_3({\textnormal{\mathversion{bold}$p$}}^T {\textnormal{\mathversion{bold}$b_3$}})),
\end{eqnarray*}
Furthermore
\begin{eqnarray*}
[p_2p_3,p_1p_3,p_1p_2]A{\textnormal{\mathversion{bold}$u$}} &=& p_2p_3+p_1p_3+p_1p_2\\
&=& 
\frac{1}{2}{\textnormal{\mathversion{bold}$p$}}^T
\left[
\begin{array}{lll}
0 & 1 & 1\\
1 & 0 & 1\\
1 & 1 & 0\\
\end{array}
\right]{\textnormal{\mathversion{bold}$p$}}\\
&=& \frac{1}{2}{\textnormal{\mathversion{bold}$p$}}^TB{\textnormal{\mathversion{bold}$p$}}.
\end{eqnarray*}
Therefore, we obtain
\begin{eqnarray*}
[p_2p_3,p_1p_3,p_1p_2]{\textnormal{\mathversion{bold}$E$}} &=& [p_2p_3,p_1p_3,p_1p_2](B{\textnormal{\mathversion{bold}$f$}}-A{\textnormal{\mathversion{bold}$u$}})\\
&=& {\textnormal{\mathversion{bold}$p$}}^T (A-\frac{1}{2}B) {\textnormal{\mathversion{bold}$p$}}.\ \mbox{\qed}
\end{eqnarray*}

Generally, the quadratic form ${\textnormal{\mathversion{bold}$p$}}^T C {\textnormal{\mathversion{bold}$p$}}$ for the matrix $C$ is equal to the symmetrical matrix $\frac{1}{2}(C+C^T)$. The above expression ${\textnormal{\mathversion{bold}$p$}}^{ T}(A-\frac{1}{2} B){\textnormal{\mathversion{bold}$p$}}$ can be written as the form with the symmetric $\frac{A+A^T-B}{2}$;  
\[
[p_2p_3,p_1p_3,p_1p_2]{\textnormal{\mathversion{bold}$E$}} = {\textnormal{\mathversion{bold}$p$}}^T\left(\frac{A+A^T-B}{2}\right){\textnormal{\mathversion{bold}$p$}}.
\]
Furthermore, by using the relation of $d_1+d_2+d_3+K+L=1$, the following result holds:
\[
[p_2p_3,p_1p_3,p_1p_2] {\textnormal{\mathversion{bold}$E$}} = {\textnormal{\mathversion{bold}$p$}}^T\left[
\begin{array}{ccc}
d_1 & \frac{d_3-d_1-d_2}{2} & \frac{d_2-d_1-d_3}{2} \\
\frac{d_3-d_1-d_2}{2} & d_2 & \frac{d_1-d_2-d_3}{2} \\
\frac{d_2-d_1-d_3}{2} & \frac{d_1-d_2-d_3}{2} & d_3
\end{array}
\right]
{\textnormal{\mathversion{bold}$p$}}
\]
The quadratic form $[p_2p_3,p_1p_3,p_1p_2]{\textnormal{\mathversion{bold}$E$}}$ is independent from the constants $K$ and $L$, and it is zero when $d_1=d_2=d_3=0$.\par
Here let us consider the case of $d_1+d_2+d_3>0$ and examine the sign of the eigenvalues of the matrix,
\[
C=\left[
\begin{array}{ccc}
d_1 & \frac{d_3-d_1-d_2}{2} & \frac{d_2-d_1-d_3}{2} \\
\frac{d_3-d_1-d_2}{2} & d_2 & \frac{d_1-d_2-d_3}{2} \\
\frac{d_2-d_1-d_3}{2} & \frac{d_1-d_2-d_3}{2} & d_3
\end{array}
\right].
\]
Since $\det C=0$ from $C{\textnormal{\mathversion{bold}$u$}}={\textnormal{\mathversion{bold}$0$}}$, the matrix $C$ has an eigenvalue $0$.
Since the trace($C$)$=d_1+d_2+d_3$ and the principal minors are
\begin{eqnarray*}
\det\left[\begin{array}{cc}
d_1 & \frac{d_3-d_1-d_2}{2}\\
\frac{d_3-d_1-d_2}{2} & d_2
\end{array}\right]
 &=& d_1d_2-\frac{(d_3-d_1-d_2)^2}{4}\\ 
 &=& \frac{2(d_1d_2+d_2d_3+d_3d_1)-(d_1^2+d_2^2+d_3^2)}{4}\\
 &=& \frac{4(d_1d_2+d_2d_3+d_3d_1)-(d_1+d_2+d_3)^2}{4},
 \end{eqnarray*}
 \begin{eqnarray*}
 \det\left[\begin{array}{cc}
d_2 & \frac{d_1-d_2-d_3}{2}\\
\frac{d_1-d_2-d_3}{2} & d_3
\end{array}\right]
 &=& \frac{4(d_1d_2+d_2d_3+d_3d_1)-(d_1+d_2+d_3)^2}{4},
 \end{eqnarray*}

  \begin{eqnarray*}
 \det\left[\begin{array}{cc}
d_1 & \frac{d_2-d_1-d_3}{2}\\
\frac{d_2-d_1-d_3}{2} & d_3
\end{array}\right]
 &=& \frac{4(d_1d_2+d_2d_3+d_3d_1)-(d_1+d_2+d_3)^2}{4},
 \end{eqnarray*}
the eigen equation of the matrix $C$ is
\begin{eqnarray*}
& &
 -t^3 + (d_1+d_2+d_3)t^2 - \frac{3}{4}\left(4(d_1d_2+d_2d_3+d_3d_1)-(d_1+d_2+d_3)^2\right)t
\\
 &=& -t\left(t^2 - (d_1+d_2+d_3)t + \frac{3}{4}\left(4(d_1d_2+d_2d_3+d_3d_1)-(d_1+d_2+d_3)^2\right)\right)
\\
&=& 0. 
\end{eqnarray*}
Thus, the other eigenvalues except $0$ are positive, when the condition
\[
H=4(d_1d_2+d_2d_3+d_3d_1)-(d_1+d_2+d_3)^2 > 0,
\]
is satisfied. From the above discussion, it can be found that in the case of $H>0$, the value of $\phi_\gamma$ in the domain ${\mathcal D}=\{{\textnormal{\mathversion{bold}$p$}} \ ; \ p_1>0,p_2>0,p_3>0 \}$ continuously decreases along the orbit if the initial value of orbit is not {\textnormal{\mathversion{bold}$p^*$}}$=[1,1,1]^T$. Furthermore, by using the Lagrange multiplier, one can check that ${\textnormal{\mathversion{bold}$p$}}^*$ gives unique only minimal value of the function $\phi_\gamma$ under the constraint $g_\gamma = 1$. The above discussion enables us to derive an important theorem on the global stability relevant to the Scarf-Hirota model.\\

\begin{theorem}
\textit{In the Scarf-Hirota model satisfying the condition (A), 
suppose $H=4(d_1d_2+d_2d_3+d_3d_1) - (d_1+d_2+d_3)^2> 0$, then the ordinary differential equations system, 
\[
\frac{dp_i}{dt}=p_i^\gamma E_i, \quad i=1,2,3\ ,
\]
is stable globally in the domain ${\mathcal D}=\{{\textnormal{\mathversion{bold}$p$}} \ ; \ p_1>0,p_2>0,p_3>0 \}$, regardless of the value of $\gamma$. Then any solution of the equations converges to the equilibrium price ${\textnormal{\mathversion{bold}$p$}}^*$.} 
\end{theorem}

\section{Local Stability}

In the differential equations system,
\[
\frac{dp_i}{dt} = p_i^\gamma E_i(p),\quad i=1,2,3\ ,
\]
the behavior of orbit near the equilibrium price ${\textnormal{\mathversion{bold}$p$}}^*$ can be determined by the eigenvalues of the Jacobian matrix, which is defined through the first-order approximation in Taylor expansion of the functions $p_i^\gamma E_i$ around ${\textnormal{\mathversion{bold}$p$}}^*$, that is
\begin{equation*}
J=\left[
\frac{\partial}{\partial p_j}(p_i^\gamma E_i)({\textnormal{\mathversion{bold}$p^*$}})
\right]_{i,j=1,2,3}.
\end{equation*}
Since the system satisfying condition(A) has the equilibrium price ${\textnormal{\mathversion{bold}$p^*$}}=[1,1,1]^T$, the elements of above matrix $J$ are 
\[
\left.\frac{\partial}{\partial p_j}(p_i^\gamma E_i)\right\vert_{{\textnormal{\mathversion{bold}$p$}}={\textnormal{\mathversion{bold}$p^*$}}}
=\left.\left(
\delta_{ij}\gamma p_i^{\gamma -1}E_i + p_i^\gamma \frac{\partial}{\partial p_j}E_i\right)\right\vert_{{\textnormal{\mathversion{bold}$p$}}={\textnormal{\mathversion{bold}$p^*$}}}=\left.\frac{\partial}{\partial p_j}E_i\right\vert_{{\textnormal{\mathversion{bold}$p$}}={\textnormal{\mathversion{bold}$p^*$}}},
\]
where $\delta_{ij}$ is Kronecker delta defined as
\[
\delta_{ij}=\left\{
\begin{array}{ll}
1 & \mbox{if } i=j\\
0 & \mbox{if } i\ne j.
\end{array}
\right.
\]
The elements of Jacobian matrix $J$ are independent of the value of $\gamma$, and they are dependent on the functions $E_i$, which are in the forms $f_1+f_2+f_3 - f_i - 1$ under the condition($A$). Since,
\begin{eqnarray*}
\frac{\partial f_h}{\partial p_j}({\textnormal{\mathversion{bold}$p^*$}}) &=&
\left.
\frac{a_{jh}({\textnormal{\mathversion{bold}$p$}}^T {\textnormal{\mathversion{bold}$b$}}_h) - ({\textnormal{\mathversion{bold}$p$}}^T {\textnormal{\mathversion{bold}$a$}}_h)b_{jh}}
{({\textnormal{\mathversion{bold}$p$}}^T {\textnormal{\mathversion{bold}$b$}}_h)^2}
\right\vert_{{\textnormal{\mathversion{bold}$p$}}={\textnormal{\mathversion{bold}$p^*$}}}\\
&=& \frac{2a_{jh}-(1-\delta_{jh})}{4},
\end{eqnarray*}
and
\[
\left.
\frac{\partial}{\partial p_j}(f_1+f_2+f_3)
\right\vert_{{\textnormal{\mathversion{bold}$p$}}={\textnormal{\mathversion{bold}$p^*$}}} = \frac{2(a_{1h}+a_{2h}+a_{3h})-(3-1)}{4} = 0,
\]
one can obtain the result
\[
\frac{\partial E_i}{\partial p_j}({\textnormal{\mathversion{bold}$p^*$}}) = -\frac{\partial f_i}{\partial p_j}({\textnormal{\mathversion{bold}$p^*$}}) = 
-\frac{2a_{ji}-(1-\delta_{ji})}{4}.
\]
Therefore, Jacobian matrix $J$ is expressed as
\[
J=\frac{1}{4}
\left[
\begin{array}{ccc}
-2a_{11} & -2a_{21}+1 & -2a_{31}+1\\
-2a_{12}+1 & -2a_{22} & -2a_{32}+1\\
-2a_{13}+1 & -2a_{23}+1 & -2a_{33}
\end{array}
\right],
\]
and clearly $J{\textnormal{\mathversion{bold}$u$}}={\textnormal{\mathversion{bold}$0$}}$ by the condition ($A$), so ${\textnormal{\mathversion{bold}$u$}}$ is an eigenvector of $J$ with the eigenvalue $0$.
From the following determinations, 
\begin{eqnarray*}
& &
\det\left[
\begin{array}{cc}
-2a_{11} & -2a_{21}+1\\
-2a_{12}+1 & -2a_{22}
\end{array}
\right]\\
 &=& 
\det \left[
\begin{array}{cc}
-2d_1 & -2(d_3+L)+1\\
-2(d_3+K)+1 & -2d_2
\end{array}
\right]\\
&=&
\det \left[
\begin{array}{cc}
-2d_1 & \left(-d_3+d_1+d_2+(K-L)\right)\\
\left(-d_3+d_1+d_2-(K-L)\right) & -2d_2
\end{array}
\right]\\
&=&
4d_1d_2 - (-d_3+d_1+d_2)^2+(K-L)^2\\
&=&-(d_1+d_2+d_3)^2+4(d_1d_2+d_2d_3+d_3d_1)+(K-L)^2 ,
\end{eqnarray*}
and
\begin{eqnarray*}
\det\left[
\begin{array}{cc}
-2a_{22} & -2a_{32}+1\\
-2a_{23}+1 & -2a_{33}
\end{array}
\right] = \det\left[
\begin{array}{cc}
-2a_{11} & -2a_{31}+1\\
-2a_{13}+1 & ^2a_{33}
\end{array}
\right]\\ = -(d_1+d_2+d_3)^2+4(d_1d_2+d_2d_3+d_3d_1)+(K-L)^2 ,\\
\end{eqnarray*}
the following eigen equation of the matrix $4J$ can be derived:
\[
-t^3 - 2(d_1+d_2+d_3)t^2 -3\hat{H}\,t
=
-t\left(t^2+2(d_1+d_2+d_3)t+3\hat{H}\right)=0,
\]
where $\hat{H}$ stands for $4(d_1d_2+d_2d_3+d_3d_1)-(d_1+d_2+d_3)^2+(K-L)^2$, which is expressed with $H$ as 
\begin{eqnarray*}
\hat{H} &=& H+(K-L)^2.
\end{eqnarray*}

The solutions of the algebraic equation
\[
t^2 +2(d_1+d_2+d_3)t+3\hat{H} = 0
\]
correspond to the other eigenvalues of $J$ except $0$. When $d_1+d_2+d_3 > 0$, the two eigenvalues satisfy the following conditions:

(ⅰ) If the discriminate is negative, that is,
\[
3\hat{H} > (d_1+d_2+d_3)^2,
\]
the two eigenvalues are complex conjugate each other while negative real part.

(ⅱ) If $0<3\hat{H} < (d_1+d_2+d_3)^2$，the two eigenvalues are negative real.

(ⅲ) If $3\hat{H} < 0$，one of eigenvalues is positive and another is negative.

The eigenvalue $0$ is not relevant to the behavior of the orbit on the surface $g_\gamma =1$, and the local stability of the orbit around ${\textnormal{\mathversion{bold}$p^*$}}$ is determined whether the condition 
\[
\hat{H} \equiv H+(K-L)^2 > 0
\]
is satisfied. The above discussion about the local stability is summarized in the following theorem.\\

\begin{theorem}
\textit{In the Scarf-Hirota model satisfying the condition (A) and  $d_1+d_2+d_3>0$, if $\hat{H}>0$, then fixed point ${\textnormal{\mathversion{bold}$p^*$}}$ is asymptotically stable in the ordinary differential equations system,}
\[
\frac{dp_i}{dt}=p_i^\gamma E_i, \quad i=1,2,3.
\]
\end{theorem}

\section{Solutions near the boundary}

In this section we study the behavior of the solutions near the boundary. 
First consider the classical case $\gamma =0$. Let $p_1 =0$, then 
\begin{eqnarray*}
E_1 &=& \frac{p_2 a_{22} + p_3 a_{32}}{p_3} + \frac{p_2 a_{23} + p_3 a_{33}}{p_2} - 1\\
&=& \frac{p_2}{p_3}a_{22} + a_{32} + a_{23} + \frac{p_3}{p_2}a_{33} -1.
\end{eqnarray*}
To find the minimum value of $E_1$ in the (open) edge ${\mathcal D}_1=\{{\textnormal{\mathversion{bold}$p$}} \ ; \ p_1=0,p_2>0,p_3>0 \}$, set
\[
\frac{p_3}{p_2} = s,
\]
and differentiate $E_1 = a_{33}s + \frac{a_{22}}{s} + a_{23}+a_{32}-1$ by $s$. Then $E_1$ takes the minimum value at $s=\sqrt{\frac{a_{22}}{a_{33}}}$, and the minimum value of $E_1$ is
\[
2\sqrt{a_{22}a_{33}} + a_{23} + a_{32} -1.
\]
Hence the necessary and sufficient condition for $E_1$ to take the positive values at the edge ${\mathcal D}_1$ is
\[
4a_{22}a_{33} > \left(1-a_{23}-a_{32}\right)^2, \mbox{　or　}
1-a_{23} - a_{32} < 0,
\]
which can be written with  $d_1,d_2,d_3, K,L$ as
\[
H>0, \mbox{ or } -d_1 + d_2 + d_3 < 0.
\]
Similarly, the necessary and sufficient conditions that $E_2$ takes a positive value at the edges ${\mathcal D}_2=\{{\textnormal{\mathversion{bold}$p$}} \ ; \ p_1>0,p_2=0,p_3>0 \}$ and $E_3$ takes a positive value at ${\mathcal D}_3=\{{\textnormal{\mathversion{bold}$p$}} \ ; \ p_1>0,p_2>0,p_3=0 \}$  are respectively
\[
H>0, \mbox{ or } d_1 - d_2 + d_3 < 0
\]
and
\[
H>0, \mbox{ or } d_1 + d_2 - d_3 < 0.
\]
Since 
\[
(-d_1+d_2+d_3)+(d_1-d_2+d_3)+(d_1+d_2-d_3) = d_1+d_2+d_3>0,
\]
the three inequalities $-d_1+d_2+d_3<0$,$d_1-d_2+d_3<0$ and $d_1+d_2-d_3<0$ can not be satisfied simultaneously, we can conclude that if $H<0$ then at least one of the following conditions is not met

(ⅰ) $E_1$ is always positive at $p_1=0$, $p_2>0$, $p_3 > 0$,

(ⅱ) $E_2$ is always positive at $p_2=0$, $p_1>0$, $p_3 > 0$, or

(ⅲ) $E_3$ is always positive at $p_3=0$, $p_1>0$, $p_2 > 0$.

 Hence if the Scarf-Hirota model with condition (A) satisfies $H<0$, then there exists a solution of the differential equation (with $\gamma =0$) which starts from the boundary of the domain $\bar{{\mathcal D}}=\{ { \textnormal{\mathversion{bold}$p$}} \ ; \ p_1 \ge 0, p_2 \ge 0 , p_3 \ge 0 \}$ and at least one price $p_i$ changes to negative, and by the continuos dependency on initial values, this implies the existence of a solution which starts from the interior of the domain $\bar{{\mathcal D}}$ and at least one price $p_i$ changes to negative~\cite{Hirota}.

Now we consider the general case $\gamma >0$. In the case $\gamma \geq 2$, the surface on which solutions are constrained by the integral does not cross the surface $p_i=0$, and in the interval $0 < \gamma <2$ there seems  to be no interesting $\gamma$ except $\gamma =1$, we study the case $\gamma =1$ only. 

Since the differential equation
\[
\frac{dE_i}{dt} = p_i E_i, \ i=1,2,3
\]
has an integral $p_1+p_2+p_3=1$, we study the solutions on the triangle 
\[
\triangle \equiv \{ {\textnormal{\mathversion{bold}$p$}} \ ; \  p_1+p_2+p_3 = 1, \ p_1\geq 0, p_2\geq 0, p_3\geq 0 \}.
\]
Since the right hand side $p_i E_i$ goes to small as $p_i$ goes close to $0$, the solutions starting from the interior of the triangle never cross the boundary. The solutions starting from the interior may only converge to the boundary.  To study that case, we first study solutions in the edge  ${\mathcal D}_1=\{{\textnormal{\mathversion{bold}$p$}} \ ; \ p_1=0,p_2>0,p_3>0 \}$. 

For simplicity, we assume that
\[
d_1>0,\ d_2>0,\ d_3>0.
\]
(Note that $d_1=0,\ d_2=0,\ d_3=0$ is in the original Scarf's model~\cite{Scarf}.)
In the edge ${\mathcal D}_1$, the Walras low $p_1E_1+p_2E_2+p_3E_3=0$ becomes $p_2E_2+p_3E_3=0$, hence
\[
E_2=0 \iff E_3=0
\]
holds.
Since the edge ${\mathcal D}_1$ is an invariant set of the differential equation, the point at $E_2=0$（hence, $E_3=0$）is a fixed point.  In the Scarf-Hirota models satisfying the condition (A), $E_i$'s satisfy\begin{eqnarray*}
E_1 &=& f_2+f_3-1\\
E_2 &=& f_3+f_1-1\\
E_3 &=& f_1+f_2-1\\
\end{eqnarray*}
and since $p_1=0$, $p_2+p_3=1$, $p_2>0$, $p_3>0$, all $f_i$ satisfy
\begin{eqnarray*}
f_1 &=& \frac{p_2 a_{21} + p_3 a_{31}}{p_2+p_3} = p_2 a_{21} + p_3 a_{31},\\
f_2 &=& \frac{p_2 a_{22} + p_3 a_{32}}{p_3},\\
f_3 &=& \frac{p_2 a_{23} + p_3 a_{33}}{p_2}.
\end{eqnarray*}
Since we assume that $d_1>0$, $d_2>0$, $d_3>0$ (hence $a_{33}>0$),
\begin{eqnarray*}
&&E_2 \to +\infty
, \mbox{when } p_1=0,\ p_2 \to 0, p_3 \to 1
\\
&&E_2 \to a_{21}+a_{23}-1 = -a_{22}>0,
\mbox{when }
p_1=0,\ p_2 \to 1,\, p_3 \to 0.
\end{eqnarray*}
These mean 
\begin{eqnarray*}
\frac{dp_2}{dt} = p_2E_2  > 0,
&\mbox{ 
when}&
p_1=0,\ p_2 \to 0, p_3 \to 1,
\\
\frac{dp_2}{dt} = p_2E_2  < 0,
&\mbox{ 
when}&
p_1=0,\ p_2 \to 1,\, p_3 \to 0,
\end{eqnarray*}
hence one can conclude  that there exists a fixed point on the edge ${\mathcal D}_1$.

Now we investigate the sign of $E_1$ at those fixed points.
Since $E_2=E_3=0$ at the fixed point, $f_2=1-f_1$,$f_3=1-f_1$, hence
\begin{eqnarray*}
E_1 &=& f_2+f_3-1\\
&=& (1-f_1)+(1-f_1)-1 = 1-2f_1.
\end{eqnarray*}
This  becomes
\[
E_1 = 1 - 2a_{31} +2(-a_{21} + a_{31})u
\]
if one set $p_2=u$, $p_3=1-u$, and also 
\[
E_1 = 1 - 2a_{21}+2(-a_{31}+a_{21})v
\]
if one set $p_3=v$, $p_2=1-v$.
The equation $f_2=f_3$, that is 
\[
p_2^2\, a_{22} + p_2p_3\,a_{32} = p_2p_3\, a_{23} + p_3^2\, a_{33},
\] 
becomes the quadratic equation
 \[
(a_{22}-a_{32}+a_{23}-a_{33})u^2 +
(a_{32}-a_{23}+2a_{33})u-a_{33} = 0,
\]
if one set $p_2=u$, $p_3=1-u$, whose leading coefficient is
\begin{eqnarray*}
a_{22}-a_{32}+a_{23}-a_{33} &=& a_{21}+a_{22}+a_{23}-a_{21}-(a_{31}+a_{32}+a_{33}-a_{31})\\
&=& 
1-a_{21}-(1-a_{31}) = -a_{21}+a_{31}.
\end{eqnarray*}
Similarly setting $p_3=v$, $p_2=1-v$, one gets the quadratic equation
\[
(a_{33}-a_{23}+a_{32}-a_{22})v^2 +
(a_{23}-a_{32}+2a_{22})v-a_{22} = 0
\]
whose leading coefficient is $-a_{31}+a_{21}$.
We set
\[p_2 = u \ \mbox{if} \ -a_{21}+a_{31} \leqq 0 \ \mbox{and} \ p_3 = v \ \mbox{if} \ -a_{31}+a_{21}<0 \]
We discuss only the case of $-a_{21}+a_{31} \leqq 0$ (the case $-a_{31}+a_{21}<0$ is similar).
Setting
$U=(-a_{21}+a_{31})u$, we rewrite the quadratic equation as
\begin{equation}\label{U}
U^2 + (a_{32}-a_{23}+2a_{33})U - (-a_{21}+a_{31})a_{33}=0
\end{equation}
and $E_1$ as
\[
E_1 = 1-2a_{31}+2U.
\]
Recall that the matrix $A$ is of the form
\[
A=\left[
\begin{array}{ccc}
d_1 & d_3+L & d_2+K\\
d_3+K & d_2 & d_1+L\\
d_2+L & d_1+K & d_3
\end{array}\right].
\]
The first order coefficient in (\ref{U}) becomes
\[
a_{32}-a_{23}+2a_{33}=2d_3+(K-L),
\]
which is positive since
\begin{eqnarray*}
0>-a_{21}+a_{31} &=& -d_3+d_2 -(K-L)\\
&=& -\left(2d_3+(K-L)\right) +d_3+d_2.
\end{eqnarray*}
The solutions of the quadratic equation are
\[
U=\frac{-(2d_3+K-L)\pm \sqrt{4d_2d_3+(K-L)^2}}{2},
\]
but only the solution satisfying
\[
0 > U \geqq -a_{21}+a_{31}
\]
gives the fixed point on the edge, since $0<u<1$, $-a_{21}+a_{31} \leqq 0$.
The inequality $U \geqq -a_{21}+a_{31}$ is equivalent to
\begin{eqnarray*}
U \geqq -a_{21}+a_{31} &\iff& -(2d_3+K-L) \pm \sqrt{4d_2d_3+(K-L)^2} \geqq 2(-a_{21}+a_{31})\\
&\iff&  -(2d_3+K-L) \pm \sqrt{4d_2d_3+(K-L)^2} \geqq 2\left(d_2-d_3-(K-L)\right)\\
&\iff& K-L-2d_2 \geqq \mp \sqrt{4d_2d_3+(K-L)^2}.
\end{eqnarray*}
Since
\[
K-L-2d_2 < \vert K-L\vert < \sqrt{4d_2d_3+(K-L)^2},
\]
the solution
\[
U=\frac{-(2d_3+K-L) - \sqrt{4d_2d_3+(K-L)^2}}{2}
\]
does not satisfy that inequality. Hence
\[
U=\frac{-(2d_3+K-L) + \sqrt{4d_2d_3+(K-L)^2}}{2}
\]
gives the required solution which gives the fixed point in the edge $\mathcal{D}_{1}$. Since at least one fixed point exists on the edge. 
Now the value of $E_1$ at that point is
\begin{eqnarray*}
E_1 &=& 1-2a_{31}+2U\\
&=& 1-2(d_2+L)-(2d_3+K-L)+ \sqrt{4d_2d_3+(K-L)^2}\\
&=&d_1-d_2-d_3 + \sqrt{4d_2d_3+(K-L)^2}
\end{eqnarray*}
Hence the condition for $E_1>0$ at that point is
\begin{itemize}
\item $d_1-d_2-d_3 > 0$（in this case $E_1>0$ for any point on the edge），or
\item $\vert d_1-d_2-d_3\vert < \sqrt{4d_2d_3+(K-L)^2}$.
\end{itemize}
Taking the square of the both sides of the later inequality, one get
\[
\hat{H}>0.
\]
Therefore, $\hat{H}>0$ implies that $E_i>0$ at the fixed point on the edge $p_i=0$ for $i=1,2,3$ and $\hat{H}<0$ does that $E_i<0$ at the fixed point on the edge, since not all the inequalities 
\begin{eqnarray*}
d_1-d_2-d_3&>&0\\
-d_1+d_2-d_3&>&0\\
-d_1-d_2+d_3&>&0.
\end{eqnarray*}
holds.
Hence we get:

\noindent
\begin{theorem}{\it
Suppose that the constants  $d_1,d_2,d_3$ in a Scarf-Hirota model with the condition {\rm (A)} satisfy
\[
d_1>0,\ d_2>0,\ d_3>0.
\]
Then the differential equation
\[
\frac{dp_i}{dt}=p_i E_i,\ i=1,2,3
\]
has locally asymptotic stable fixed points on the edges if $\hat{H} < 0$,
and it has no locally asymptotic stable points on the edges if $\hat{H} > 0$
}

\end{theorem}

From the Theorem 2, the equilibrium is globally asymptotic stable if $H>0$, and from the Theorem 3, it is not even locally asymptotic stable if $\hat{H}<0$. The remaining case is $H<0$, $\hat{H}>0$. In the case of $\gamma =0$, there emerge solutions starting from the inside and going to the outside as soon as the condition $H>0$ breaks, but in the case $\gamma=1$, there is no solutions which go to the outside, and from Theorem 4 no solution converges to a point on the edge as long as $\hat{H}>0$ holds. Therefore if there exist parameters $d_1,\,d_2,\,d_3,\,K,\,L$ for which the global asymptotic stability condition does not holds and the inequalities $H<0$ and $\hat{H}>0$ hold, then there exists a locally asymptotically stable periodic orbit for that parameters by the Poincar\'{e}-Bendixson Theorem.

But at least there seems to be no reason that the inequality $H<0$ assures the appearance of a periodic orbit, so it is natural to expect that there exists a globally asymptotic stable case with $H<0$, $\hat{H}>0$. In fact, the following example of Figure.~1 is of that case. In this example where the matrix A is set by the parameters, $d_1=0.1$, $d_2=0.1$, $d_3=0.5$, $K=0$ and $L=0.3$,
\begin{figure}[htbp]
 \begin{center}
    \includegraphics[keepaspectratio=true,height=70mm]{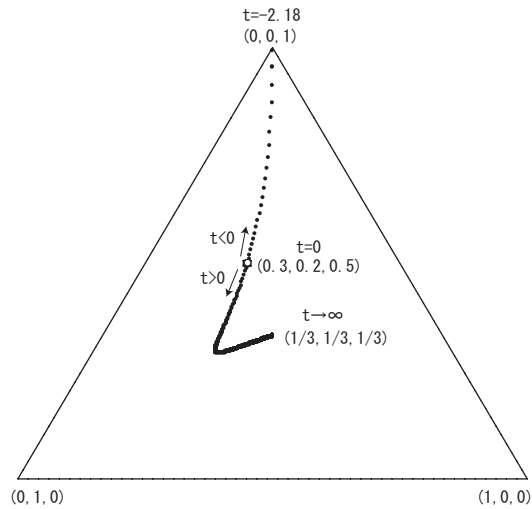}
  \end{center}
  \caption{The solution that starts from the initial value $(0.3,0.2,0.5)$. }
 \label{fig:FIGUREY.eps}
\end{figure}
the solution starting from the initial value $(0.3,0.2,0.5)$ converges to the equilibrium $(\frac{1}{3},\frac{1}{3},\frac{1}{3})$ as $t\to \infty$, and as time goes to negative, the solution reaches to the node $(0,0,1)$ at $t=-2.18$ (in the case of $\gamma \leq 1$, the right hand side of differential equation is not defined at the nodes, which makes it possible for the solution to reach a node before $t\to -\infty$). This example seems only to say that a solution starting from a particular initial value converges to the equilibrium, However the existence of such a solution is sufficient to conclude that all solutions starting from an arbitrary initial value inside the triangle converges to the equilibrium. If there exists a solution which does not converge to the equilibrium, then by the Poincare-Bendixson Theorem, the solution winds toward a periodic orbit. By the fixed point theorem, there exists a equilibrium in the domain bounded by the periodic orbit. However $p*$ is the unique equilibrium in the interior, the periodic orbit separates $p*$ from the edges of the triangle; this is impossible since the above example gives a solution which connects one of the edge and the equilibrium.

We have tested for many parameters with $H<0$, $\hat{H}>0$, but we can not find any example which is not global asymptotically stable, and we conjecture: "local asymptotically stability implies global asymptotically stability in the price-scaled price adjustment process".

Fig1: The solution that starts from the initial value $(0.3,0.2,0.5)$.
\end{document}